\documentclass{article}
\usepackage{amsmath}
\usepackage{amssymb}
\usepackage{amsthm}
\usepackage[margin=1in]{geometry}
\usepackage[utf8]{inputenc}
\usepackage{graphicx}
\usepackage{caption}
\usepackage{subcaption}
\usepackage{float}
\usepackage{setspace}
\usepackage[titletoc,title]{appendix}
\usepackage{xcolor}
\allowdisplaybreaks
\usepackage[noblocks]{authblk}
\usepackage{amsthm}
\usepackage{tikz}
\usetikzlibrary{shapes,arrows,chains}
\usepackage{mathtools}
\usepackage{soul}
\usepackage[style=authoryear]{biblatex}
\usepackage{blkarray}
\usepackage{pgfplots}

\title{Synchronization in a Two-supplier Assembly System: Combining a Fixed Lead-time Module with Capacitated Make-to-Order Production}
\author{M.S. Meijer \thanks{Corresponding author. Email: m.s.meijer@tue.nl}}
\author{W. van Jaarsveld}
\author{A.G. de Kok}
\affil{Department of Industrial Engineering and Innovation Sciences, Eindhoven University of Technology, Eindhoven, The Netherlands}
\date{April 2021}

\newtheorem{theorem}{Theorem}
\newtheorem{observation}{Observation}
\newtheorem{lemma}{Lemma}

\newcommand{\compone}{$m_1$}
\newcommand{\comptwo}{$m_2$}
\newcommand{\waitingcost}{$b$}
\newcommand{\numberevents}{$X$}
\newcommand{\startvector}{$s$}
\newcommand{\production}[1]{$M_1(#1 )$}
\newcommand{\transitionmatrix}{$\bar{P}$}

\begin{document}

\maketitle

\begin{abstract}
A high-tech manufacturer often produces products that consist of many modules. These modules are either sourced from one of its suppliers or produced in-house. In this paper we study an assembly system in which one module is sourced from a supplier with a fixed lead-time, while the other module is produced by the manufacturer itself in a make-to-order production system. Since unavailability of one of the modules has costly consequences for the production of the end-product, it is important to coordinate between the ordering policy for one module and the production of the other. We propose an order policy for the lead-time module with base-stock levels depending on the number of outstanding orders in the production system of the in-house produced module. We prove monotonicity properties of this policy and show optimality. Furthermore, we conduct a computational experiment to evaluate how the costs of this policy compare to those of a policy with fixed base-stock levels and show that average savings of 17\% are attained.  
\end{abstract}

\section{Introduction}

High-tech original equipment manufacturers (OEMs) produce complex products composed of many different modules that are either produced by the OEM itself or sourced from one of its suppliers. To be able to assemble the final product and deliver it to the customer, the OEM needs to organize its production and ordering activities such that all modules are available at the point of assembly and inventory holding and waiting time costs are minimized. The main challenge herein is that external suppliers deliver after an agreed lead time and the OEM needs to align its production with this.
Since many modules have long lead times and uncertainties are present also in the in-house production process, this is a complicated problem.

In this paper, we consider a small-scale assembly system with an end-product consisting of two modules. One of the modules is produced in-house by the OEM after a customer order has been placed (make-to-order). The second module is sourced from a supplier with a given lead time. Customers want to receive the product as soon as possible after placing their order. Therefore, the OEM wants to start assembly as soon as the in-house produced component is ready, to avoid costly customer waiting times. However, both modules need to be available at the time of assembly. If production of the first module has finished while the second module is not yet available, the first module needs to be stored and the customer has to wait. On the other hand, if the second module is on stock but the in-house produced module is not ready, the second module needs to be stored. 
To avoid high inventory holding and waiting time costs, the OEM needs to synchronize the production of the make-to-order component with the order policy of the component sourced from its supplier.

There are different approaches to modeling supply processes. Lead-times can be deterministic and thus independent of the volume of an order. This is for example the case for bulk products that are shipped overseas with a given transport time. For some products, however, the lead-time is highly dependent on the volume. Important work on coordination of ordering decisions for multiple items with a deterministic lead-time was done by \textcite{rosling1989}. Coordination of multiple items with stochastic lead-times was studied by among others \textcite{benjaafar2006}. However, literature is lacking on coordination between items with deterministic and stochastic lead-times, which is the focus of this paper. 
There are many practical examples of assembly systems in which this type of coordination is required. The production of lithography machines is a capacitated process for which modules with a long lead-time, such as lenses, are sourced at external suppliers. Since these lenses have a long, but predictable, lead-time, it is important to synchronize the ordering of these modules with the capacitated production process. 
More generally, final assembly is often a capacitated process and thus has a stochastic lead-time, after which some final items from external sources with a fixed lead-time need to be added, such as packaging materials. Accessories or options packages are also often added after the main product has been produced.

We consider two different models of an assembly system: a continuous time model and a discrete time model. In the continuous time model, we assume Poisson distributed demand and a single-server production system with exponential production times for the in-house produced module.
Once a customer order arrives, production of the module can start as soon as there is available capacity. Until capacity becomes available, the order has to wait in the queue. We aim to synchronize the output process of this production system with the order policy of the other module, to avoid large inventories that give rise to holding costs or penalty costs incurred for waiting customers. Useful information on the expected production of the MTO module, and thus the required units of the second module, can be obtained from the number of customer orders waiting for production of the MTO module. 
Since the production capacity of the MTO module is fixed and assembly starts as soon as both modules are available, the inventory levels of both finished modules are influenced by the inventory position of the module sourced from the supplier.
We propose a base-stock policy where the target inventory position of the supplier-sourced module is dependent on the number of customer orders waiting for production of the MTO module. We prove monotonicity properties of this policy and show optimality.
Numerical results demonstrate that the proposed policy can generate considerable savings compared to a base-stock policy with fixed base-stock levels.
Next, we investigate whether the results we obtained for the continuous time model can be extended to discrete time. We consider a discrete time model with a production capacity per period that can either be random or fixed. This model would be more suitable when there is a fixed number of products that can be produced per period or when available equipment has a random yield per period. 
Also in this setting, we can prove optimality of the base-stock policy with target inventory positions depending on the number of outstanding orders. 
Numerical results again indicate that considerable savings can be attained by considering information on the number of waiting orders for the MTO module. However, we also observe some differences in the results compared to the continuous time model. For example, in continuous time we observed an increase in the average percentage savings as lead-time of increased, whereas in discrete time we observe a decrease.

This paper is organized as follows. We review relevant literature in \S\ref{sec: literature}. In \S\ref{sec: model} we explain our model in detail. Using this model we derive the optimal base-stock policy for the module sourced from the supplier and show how to compute the policy parameters based on the state of the in-house production system. A computational experiment is provided in \S\ref{sec: computational analysis}, which shows an example of what the inventory policy will look like and compares the expected costs of this state-dependent policy to those of a policy with a fixed base-stock level. In \S\ref{sec: discrete model}, we show optimality of this policy also for the discrete time model. The paper is concluded in \S\ref{sec: conclusions}.

\section{Literature Review}\label{sec: literature}

Assembly systems have been studied extensively, for example by \textcite{schmidt1985} who provide optimal policies for assembly systems with two components.
\textcite{kok2018} review the extensive literature on multi-echelon inventory management over the past decades, covering convergent, divergent and more general structures and any combination of make-to-order, make-to-stock and assemble-to-order production.
\textcite{atan2017} provide a recent overview of literature studying assemble-to-order systems. They state that the main challenge in continuous review models with a single end-product is the synchronization of component orders. 
We can classify literature on assembly systems in two groups with respect to capacity constraints. Literature in the first group does not take into account capacity constraints. They assume stochasticity on the demand side, but deterministic lead times and unlimited supply. 
In the second group we find literature concerning assembly systems with capacity constraints, resulting in stochastic lead times.

First, we will focus on literature that studies coordination in uncapacitated assembly systems.
\textcite{rosling1989} shows that under certain conditions a multi-stage assembly system with fixed assembly times can be reduced to an equivalent serial system, for which optimal policies are designed by \textcite{clark1960}. 
More recently, variations of assembly systems have been studied, such as systems with components that have both different lead times and review periods \parencite{karaarslan2018}. 
\textcite{martinez2010} provide a closed-form formula to determine whether or not an order should be placed. Additionally, they provide conditions for optimality of such myopic policy. 
\textcite{lu2015} derive (asymptotically) optimal policies for both inventory replenishment and inventory allocation for assemble-to-order $N$- and $W$-systems.

Capacity constraints are often modeled using finite-capacity queueing systems.
\textcite{song1999} consider a production system with exponential production times in single-server queues with a finite queue and show how to obtain performance measures.
The past decades, such assembly systems with finite production capacities are studied increasingly \parencite[e.g.][]{bollapragada2015,elhafsi2010,plambeck2008,toktacs2011}. \textcite{benjaafar2006} study an assembly system consisting of $m$ components required to satisfy demand of $n$ customer classes. A policy needs to specify when to produce each component and whether or not to satisfy incoming customer orders from on-hand inventory. They show a base-stock policy with dynamic base-stock levels is optimal. \textcite{benjaafar2011} extend this work to the case where production facilities do not only produce components, but also sub-assemblies. 
\textcite{huh2010} focus on base-stock policies, as these are often used in practice. They show convexity of the shortage costs with respect to the order-up-to levels and discuss algorithmic implications.
\textcite{cheng2011} study a problem with unpredictable machine breakdowns and endogenous load-dependent lead-times.
\textcite{song1993} study an assembly system with stochastic lead-times and Markov-modulated demand, meaning that demand rates are dependent on the state of an underlying variable. 
Several variations and extensions of the work of \textcite{song1993} have been considered, including the work of \textcite{chen2001}. \textcite{gallego2004} consider, besides a Markov-modulated demand process, also a Markov-modulated supply process that was driven by an independent Markov chain. Furthermore, they consider finite production capacities.
Similarly, \textcite{mohebbi2006} studies a situation where supply and demand are subject to independent random environmental conditions where production up to the storage capacity is initiated as soon as the inventory level drops below the limit. 
\textcite{muharremoglu2008} study a serial system with multiple stages and stochastic lead times with Markov-modulated demand. They provide an approach for decomposing the serial inventory problem into decoupled subproblems each consisting of a single unit and a single customer. They show that state-dependent base-stock policies are optimal and provide an efficient algorithm to compute the base-stock levels. 

Our work combines the two research streams discussed above by studying coordination in an assembly system that combines a module sourced from an unconstrained supplier with fixed lead-time with a module that is produced in a capacitated system with stochastic lead-time. Due to the prevalence of combinations of two such supply streams in practice, this is a relevant topic to study. However, besides its practical relevance, this problem is interesting from a theoretical perspective. Clearly, theory on uncapacitated systems with deterministic lead-times cannot be generalized to situations in which one supply stream has uncertain lead-time, since it is no longer possible to order items based on their lead-time. Furthermore, contrary to most studies on assembly with stochastic lead-times, we assume the capacity investment decision for the stochastic production system to be fixed and coordinate the availability of both modules through the inventory policy of the supplier-sourced module.

\section{Continuous Time Model}\label{sec: model}
In this section we consider a continuous time model of the assembly system. In \S \ref{sec: discrete model} we will consider the case where the assembly process is modeled in discrete time.

Consider a high-tech end product that is composed of two modules. The first module, denoted by \compone, is customer-specific and made to order by the OEM itself. Production of \compone{} starts as soon as capacity is available after arrival of the customer order.
The second module, denoted by \comptwo, is sourced from a supplier with lead-time $L$. 
Demand of the end product is modeled as a Poisson process with customers arriving at rate $\lambda$. We assume that \compone{} is produced in a single-server queue with exponential service times with rate $\mu$. 
As a consequence, the production of module \compone{} evolves as an $M/M/1$ queue. A departure from the queuing system then represents a finished \compone{} module that can be used to assemble the end product. 

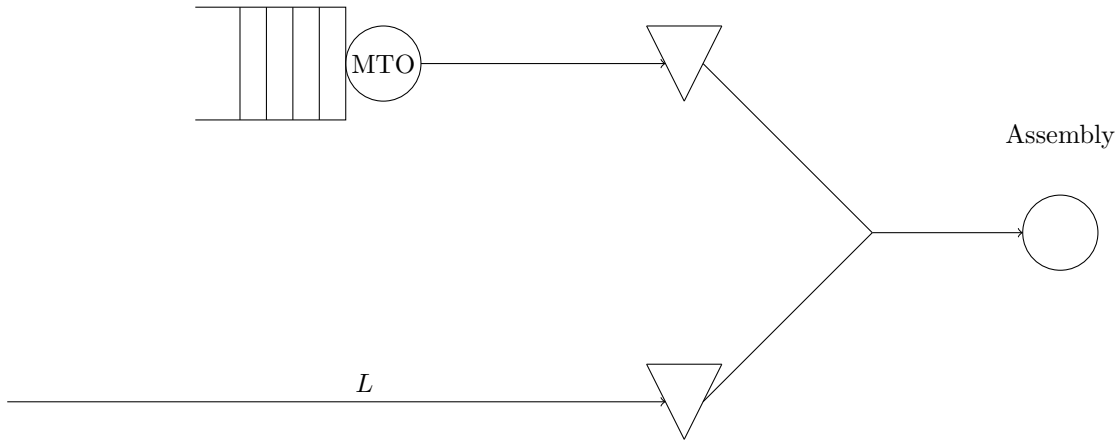
\begin{figure}[ht]
\centering
\scalebox{1}{
\begin{tikzpicture}
    \draw (5cm,2.25cm) circle [radius=0.5cm] node{MTO};
    \draw (2.5,3) -- ++(2cm,0) -- ++(0,-1.5cm) -- ++(-2cm,0);
    \foreach \i in {1,...,4}
        \draw (4.5cm-\i*10pt,3) -- +(0,-1.5cm);
    
    \draw[->](5.5,2.25)--(8.75,2.25);
    
    \draw(8.5,2.75)--(9.5,2.75);
    \draw(8.5,2.75)--(9,1.75);
    \draw(9.5,2.75)--(9,1.75);
    
    \draw[->](0,-2.25)--(8.75,-2.25)node[above,yshift=0cm,xshift=-4cm]{\text{$L$}};
    
    \draw(8.5,-1.75)--(9.5,-1.75);
    \draw(8.5,-1.75)--(9,-2.75);
    \draw(9.5,-1.75)--(9,-2.75);
    
    \draw(9.25,2.25)--(11.5,0);
    \draw(9.25,-2.25)--(11.5,0);
    \draw[->](11.5,0)--(13.5,0);
    \draw (14,0) circle [radius=0.5cm] node[above,yshift=1cm,xshift=0cm]{\text{Assembly}};
\end{tikzpicture}}
\caption{Sketch of the assembly system}
\label{fig: sketch assembly system}
\end{figure}

When module \compone{} is finished and module \comptwo{} is available, the final product can be assembled. If module \comptwo{} arrives before \compone{} is available, it needs to be stored and holding costs are incurred. If module \comptwo{} has not yet arrived when module \compone{} is finished, \compone{} needs to be stored. 
A sketch of this system is given in Figure \ref{fig: sketch assembly system}.
To formulate our model of the given assembly system, we introduce the additional notation given in Table \ref{tab: notation}.

\begin{table}[ht]
    \centering
    \caption{Notation}
    \begin{tabular}{ll}
    \hline
        $I_1^P(t)$ & orders of \compone{} being processed or waiting to be processed at time $t$\\
        $I_1(t)$ & inventory of finished \compone{} at time $t$ ready for assembly\\
        $I_2^T(t)$ & ordered \comptwo{} modules that are in transit at time $t$\\
        $I_2(t)$ & inventory of \comptwo{} at time $t$\\
        $I^A(t)$ & final products being assembled at time $t$\\
        $IP_2(t)$ & inventory position of module \comptwo{} at time $t$\\
        \production{t} & cumulative production of module \compone{} until time $t$\\
        $h_1$ & unit holding costs for \compone\\
        $h_2$ & unit holding costs for \comptwo\\
        $b$ & costs for customers waiting for their final product\\
        $L$ &  lead-time for \comptwo\\
        \hline
    \end{tabular}
    \label{tab: notation}
\end{table}

Customers want to receive the final product as soon as possible after placing their order, which is represented by waiting costs \waitingcost. Since module \compone{} is made to order and therefore there is a one to one correspondance between customer orders and \compone{} modules, we can formulate the cost function at time $t$ consisting of four parts:
\begin{enumerate}
    \item \compone{} orders that are waiting to be processed represent customers that are waiting, hence cost \waitingcost{} is incurred for every item.
    \item Similarly, \compone{} modules that are finished, waiting to be assembled also represent customers that are waiting. Additionally, holding costs are incurred for the finished modules, leading to costs $b+h_1$ per unit.
    \item \comptwo{} modules that are delivered and are waiting to be merged with module \compone{} give rise to holding costs $h_2$.
    \item Final products that are being assembled consist of \compone{} and \comptwo{}. Also, every final product is coupled to a customer order and thus represents a waiting customer. Therefore, costs $b+h_1+h_2$ are incurred.
\end{enumerate}
Combining these parts gives cost function
\begin{equation}\label{eq: costs}
    C(t) = b I_1^P(t) + (b+h_1)I_1(t) + h_2 I_2(t) + (b+h_1+h_2) I^A(t).
\end{equation}
Since the capacity for producing \compone{} and for assembly are fixed, $I_1^P(t)$ and $I^A(t)$ are considered to be exogenous as $t\rightarrow\infty$ so that  
we can write $E[C(t)] = E[\tilde{C}(t)]+[constant]$, with 
\begin{equation}
    \tilde{C}(t) = (b+h_1)I_1(t) + h_2 I_2(t)
\end{equation}
being the cost function of interest. We are interested in minimizing the average costs of operating the assembly system over time. 

\subsection{Inventory policy for module \comptwo{}}\label{sec: inventory policy c2}
Since minimizing the average costs over time is not straightforward, we will first try to solve a related problem of minimizing the expected costs at a certain point in time. After obtaining the solution that minimizes the expected costs at a fixed point in time, we will return to the problem of minimizing the average costs over time. Since decisions at time $t$ affect costs at time $t+L$, we consider the expected costs $\tilde{C}(t+L)$ given the relevant information about the state of the system available at time $t$:
\begin{equation}\label{eq: expected costs}
    E[\tilde{C}(t+L)|I_1(t), I_2(t), I_2^T(t), I_1^P(t)]=E[\tilde{C}(t+L)|I_1(t), I_2(t)+ I_2^T(t), I_1^P(t)].
\end{equation}
The equality holds because all \comptwo{} in transit at time $t$ will have been delivered at time $t+L$. Since $\tilde{C}(t+L)$ is determined by the values of $I_1(t+L)$ and $I_2(t+L)$, we are interested in the values of $I_1(t+L)$ and $I_2(t+L)$ given the information available at time $t$ ($I_1(t)$, $I_2(t)+ I_2^T(t)$ and $I_1^P(t)$). 

Since \production{t} denotes the cumulative production of module \compone{} until time $t$, 
it follows that the production during time interval $[t,t+L]$ can be written as \production{t+L}$-$\production{t}. Since module \compone{} and \comptwo{} are combined into the final product, the inventory level of module \compone{} can decrease by at most the availability of module \comptwo, which is equal to $I_2(t)+I_2^T(t)$. This means that we can write the inventory of module \compone{} at time $t+L$ as
\begin{equation*}
    I_1(t+L) = \max \left\{I_1(t) + \textrm{\production{t+L}} - \textrm{\production{t}} - I_2(t) - I_2^T(t),0\right\}.
\end{equation*}
Similarly, we can write
\begin{equation*}
    I_2(t+L) = \max \left\{I_2(t) + I_2^T(t) - I_1(t) - \left(\textrm{\production{t+L}}- \textrm{\production{t}}\right),0\right\}.
\end{equation*}

The inventory position of module \comptwo{} increases when a new order is placed. When production of a module \compone{} is finished, the module is ready to be merged with module \comptwo{}, leading to a decrease in the inventory position of \comptwo. Therefore, the inventory position of module \comptwo{} at time $t$ is equal to $IP_2(t)=I_2(t) + I_2^T(t) - I_1(t)$. This allows us to rewrite the expression for expected costs given in Equation \eqref{eq: expected costs} as:
\begin{align}\label{eq: expected costs IP}
    E[\tilde{C}(t+L)|IP_2(t), I_1^P(t)]=E[(b+h_1)(\textrm{\production{t+L}} - \textrm{\production{t}} - IP_2(t) )^+ + h_2 (IP_2(t) - (\textrm{\production{t+L}} - \textrm{\production{t}}))^+|I_1^P(t)].
\end{align}
The expected cost is thus determined by the production of module \compone{} and the inventory position of module \comptwo, which is controlled by the order policy. 

Since the capacity for producing module \compone{} is fixed and hence the production of \compone{} cannot be influenced, the inventory control policy for module \comptwo{} is the only decision that can influence the inventory levels of both modules and thus the expected costs. 
Therefore, we aim to find the inventory policy for module \comptwo{} that minimizes the expected costs as given in Equation \eqref{eq: expected costs IP}.
We will consider a myopic inventory policy, where at any time $t$ the target inventory position of module \comptwo{} is determined that minimizes $E[\tilde{C}(t+L)|IP_2(t), I_1^P(t)]$. This yields the minimization problem given in Equation \eqref{eq: minimization}.
\begin{align}\label{eq: minimization}
    \min_{IP_2(t)} E[(b+h_1)(\textrm{\production{t+L}} - \textrm{\production{t}} - IP_2(t) )^+ + h_2 (IP_2(t) - (\textrm{\production{t+L}} - \textrm{\production{t}}))^+|I_1^P(t)].
\end{align}
We denote the minimizing target inventory position by $\tilde{IP}_2(t)$. Observation \ref{obs: target IP function of I_1^P(t)} allows us to write $\tilde{IP}_2(t)=\tilde{IP}_2(I_1^P(t))$. This myopic inventory policy thus consists of a list of target inventory position levels for every current state of the system.

\begin{observation}\label{obs: target IP function of I_1^P(t)}
The target inventory position of module \comptwo{} is a function of the number of waiting orders in the queue for \compone.
\end{observation}

Next, we want to assess whether this myopic policy is a good policy. For this purpose, we will further analyze how the target inventory position $\tilde{IP}_2(I_1^P(t))$ responds to changes in $I_1^P(t)$.
First of all, in Theorem \ref{thm: monotonicity} we show that the target inventory position of module \comptwo{} is non-decreasing in the number of customers in the queue. If this would not hold, there could be situations in which the target inventory position is lower than the current inventory position and one would like to place a negative order. Similarly, Theorem \ref{thm: max_increase} shows that the increase of $\tilde{IP}_2(I_1^P(t))$ is at most one when $I_1^P(t)$ increases by one. All proofs are given in the Appendix.

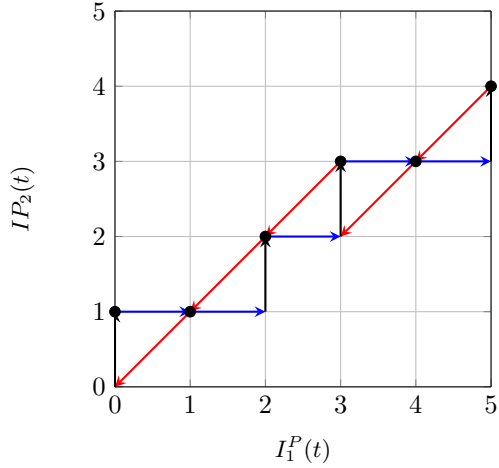
\begin{figure}
    \centering
\begin{tikzpicture}
\centering
\begin{axis}[%
footnotesize,
width=5cm,
height=5cm,
axis equal,
scale only axis,
separate axis lines,
grid=both,
every outer x axis line/.append style={black},
every x tick label/.append style={font=\color{black}},
xmin=0,
xmax=5,
xlabel={$I^P_1(t)$},
every outer y axis line/.append style={black},
every y tick label/.append style={font=\color{black}},
ymin=0,
ymax=5,
ylabel={$IP_2(t)$},
minor tick num=0,
xtick distance=1,
ytick distance=1,
axis background/.style={fill=white},
legend style={at={(0.97,0.03)},anchor=south east,legend cell align=left,align=left,fill=white}
]
\addplot[only marks, color=black, mark=*, mark options={scale=1}] coordinates {
    (0,1) (1,1) (2,2) (3,3) (4,3) (5,4)};
\addplot[blue,thick,
        quiver={u=\thisrow{u},v=\thisrow{v}},
        -stealth] 
	table 
	{
	x y u v
	0 1 1 0
	1 1 1 0
	2 2 1 0
	3 3 1 0
	4 3 1 0
	};
\addplot[red,thick,
        quiver={u=\thisrow{u},v=\thisrow{v}},
        -stealth] 
	table 
	{
	x y u v
	1 1 -1 -1
	2 2 -1 -1
	3 3 -1 -1
	4 3 -1 -1
	5 4 -1 -1
	};
\addplot[black,thick,
        quiver={u=\thisrow{u},v=\thisrow{v}},
        -stealth] 
	table 
	{
	x y u v
	0 0 0 1
	2 1 0 1
	3 2 0 1
	5 3 0 1
	};
\end{axis}
\end{tikzpicture}
\caption{Illustration of order policy. Blue arrows correspond to customer order arrivals, black arrows to orders of \comptwo{} and red arrows to assembly.}
\label{fig: illustration policy}
\end{figure}

\begin{theorem}\label{thm: monotonicity}
$\tilde{IP}_2(I_1^P(t))$ is monotonically non-decreasing in $I_1^P(t)$.
\end{theorem}
\begin{theorem}\label{thm: max_increase}
If an additional customer enters the system, the target inventory position of module \comptwo{} increases by at most 1, i.e. $\tilde{IP}_2(I_1^P(t)+1)-\tilde{IP}_2(I_1^P(t))\leq1$
\end{theorem}
Combining Theorem \ref{thm: monotonicity} and Theorem \ref{thm: max_increase}, we conclude that if $I_1^P(t)$ increases by one, the target IP level remains the same or increases by one. The proposed myopic inventory policy is thus a well-defined policy. 
We give an illustration of this in Figure \ref{fig: illustration policy}. On the horizontal axis we have the inventory position of module \compone{} and on the vertical axis that of module \comptwo. The black dots indicate the target inventory position of \comptwo{} for the given $I^P_1(t)$. The red arrows correspond to an assembly step where both \compone{} and \comptwo{} are used, reducing both $I_1^P(t)$ and $IP_2(t)$ by one. The blue arrows correspond to a customer order arrival, leading to an increase in $I^P_1(t)$ of 1. The black arrows show the required orders of module \comptwo{} to reach the target inventory position. We observe that both the blue and the red arrows never go to a point where the inventory position of module \comptwo{} exceeds it target. Consequently, all black arrows correspond to an order of size 1. Hence, this policy behaves well in all situations, in the sense that it never prescribes negative orders.

From Theorems \ref{thm: monotonicity} and \ref{thm: max_increase} we learn that the proposed policy that minimizes costs at a specific time point is well-defined. The next question is how this policy performs when we return to the original problem in which we aim to minimize the average costs over time.
When we consider the optimal inventory policy that minimizes the average costs there also is a corresponding inventory position $IP_2^*(t)$ at every time $t$.
By definition, the myopically optimal target inventory position minimizes the costs at every point in time. 
Therefore, we can conclude that this myopic policy is also optimal for minimizing average costs over time and thus that $IP_2^*(t)=\tilde{IP}_2(I_1^P(t))$. This is formalized in Theorem \ref{thm: optimality}. 

\begin{theorem}\label{thm: optimality}
The myopic inventory policy for module \comptwo{} is optimal for minimizing average costs.
\end{theorem}

\subsection{Computing $\tilde{IP}_2(I_1^P(t))$}
Now that we have characterized the myopic inventory policy, we are interested in how the target inventory level $\tilde{IP}_2(I_1^P(t))$ can be computed for a given value of $I_1^P(t)$. This requires us to take a closer look at the rate at which module \comptwo{} is used in the assembly process and the cost structure. 

When module \compone{} is finished and module \comptwo{} is available, the final product can be assembled. However, if module \comptwo{} arrives before \compone{} is available it needs to be stored and holding costs are incurred and if module \comptwo{} has not yet arrived when module \compone{} is finished, \compone{} needs to be stored and the customer is waiting.
In other words, the minimization problem given in Equation \eqref{eq: minimization} is a Newsvendor problem with shortage costs for module \comptwo{} equal to $b+h_1$ and  overage costs $h_2$.

The requirement for module \comptwo{} in the assembly system during period $[t,t+L]$ is equal to the production of module \compone{} during $[t,t+L]$, which is denoted by $\textrm{\production{t+L}} - \textrm{\production{t}}$. Since the production during $[t,t+L]$ depends on the number of waiting customers at time $t$, the distribution of $\textrm{\production{t+L}} - \textrm{\production{t}}$ depends on $I_1^P(t)$. Therefore, when synchronizing the order policy of  module \comptwo{} with the production of module \compone, we need to consider the distribution of $\textrm{\production{t+L}} - \textrm{\production{t}}$ given $I_1^P(t)$.
We then want to determine $\tilde{IP}_2(I_1^P(t))$ such that $P(\textrm{\production{t+L}} - \textrm{\production{t}}<\tilde{IP}_2(I_1^P(t)))\geq CF$, where $CF=\frac{b+h_1}{b+h_1+h_2}$ is the critical fractile. 

To analyze the production of \compone, we model the state of the inventory system, consisting of the number of waiting customer orders and the number of finished modules \compone, as a Markov process. Let $(i,j)$ denote the state of the system with $i$ the number of \compone{} jobs waiting or currently in process in the system and $j$ the number of \compone{} modules produced, with $i=0,1,2,\ldots$ and $j=0,1,2,\ldots$. 
A customer order arrival occurs with probability $\tilde{\lambda}=\frac{\lambda}{\lambda+\mu}$ and generates a transition from $(i,j)$ to $(i+1,j)$. An exit corresponds to finished production of a module \compone{}. This occurs with probability $\tilde{\mu}=\frac{\mu}{\lambda+\mu}$. When there are customers in the system ($i>0$), a completion of module \compone{} means that the process moves to state $(i-1,j+1)$. When $i=0$, meaning that there are no customer orders in the system, we do allow for exits, but the system then remains in the same state.
This yields the following transition probabilities:
\begin{align*}
    P_{(i,j),(i+1,j)}=\tilde{\lambda}, &\\
    P_{(i,j),(i,j)}=\tilde{\mu} & \textrm{ if } i=0, \textrm{ and}\\
    P_{(i,j),(i-1,j+1)}=\tilde{\mu} & \textrm{ if } i>0.
\end{align*}
We are interested in the number of modules \compone{} produced by time $t+L$, given that there are $I_1^P(t)$ jobs waiting in the system at time $t$. In other words, we are interested in the value of $j$ after $L$ time units, starting from state $(I_1^P(t),0)$.

In order to analyze the production during $[t,t+L]$, we condition on the total number of transitions, consisting of both customer arrivals and finished production of a module, during time period $[t,t+L]$, denoted by \numberevents{}.
When we condition on \numberevents, both the number of customers in the queue and the number of modules \compone{} that are produced during time interval $[t,t+L]$ are bounded. As a result the state space is bounded and we can define the transition probability matrix, which we denote by \transitionmatrix. Let \startvector{} be a vector with length equal to the number of states with all zero entries, except for a one at the state corresponding to $(I_1^P(t),0)$. The probability distribution over the state space after \numberevents{}$=k$ events is then given by \startvector$\text{\transitionmatrix}^k$.
Since the production of \compone{} is modeled as an $M/M/1$ queue, \numberevents{} has a Poisson distribution with parameter $(\lambda+\mu)L$, i.e. 
\begin{equation*}
    P(X=k)=\frac{((\lambda+\mu)L)^k}{k!}e^{-(\lambda+\mu)L}.
\end{equation*}
For each realization of $X=k$ we can determine the probability distribution over the output of \compone{} from $\text{\startvector}\bar{P}^k$, from which we obtain the unconditional distribution over the state space after $L$ time units starting from state $(I_1^P(t),0)$. 
Since the support of the Poisson distribution is all natural numbers starting from 0, i.e. $k\in\mathbb{N}_0$, we need to bound the states space in order to obtain closed-form expressions.
Therefore, we introduce an upper bound on the number of transitions during $[t,t+L]$, denoted by \numberevents$^U$, such that $P($\numberevents$\geq$\numberevents$^U)$ is negligible.
Then the distribution over the number of modules \compone{} produced in the time interval $[t,t+L]$ can be obtained from Equation \eqref{eq: unconditional_distribution}.
\begin{equation}\label{eq: unconditional_distribution}
    P(\textrm{\production{t+L}} - \textrm{\production{t}}=j|I_1^P(t)=i)=\sum_{k=0}^{\text{\numberevents}^U} P(X=k)P(\textrm{\production{t+L}} - \textrm{\production{t}}=j|X=k,I_1^P(t)=i)
\end{equation}
We specify \numberevents$^U$ using Cantelli's inequalities as described by \textcite{ghosh2002}. Since \numberevents$\sim Poisson((\lambda+\mu)L)$, this gives the upper bound provided in Lemma \ref{lem: upper bound}.
\begin{lemma}\label{lem: upper bound}
For any $0<\epsilon<1$, $P(X\ge X^U)\leq \epsilon$ holds if $X^U=(\lambda+\mu)L+\sqrt{\left(1-\frac{1}{\epsilon}\right)(\lambda+\mu)L}$.
\end{lemma}

To illustrate this procedure, we will now provide a small scale example that shows the intuition behind this approach.


\subsubsection{Small scale example production \compone}
We consider a small example in which the parameters are such that the probabilities of having more than 2 customers in the queue or of producing more than 2 units of module \compone{} during $[t,t+L]$ are negligible.
This gives the following states and  transition probabilities, where $(0,3)$ is added as an absorbing state such that the probabilities in each row add up to one: 

\begin{equation*}
\bar{P}=   \begin{blockarray}{ccccccccccc}
&(0,0) & (1,0) & (2,0) & (0,1) & (1,1) & (2,1) & (0,2) & (1,2) & (2,2) & (0,3) \\
\begin{block}{c[cccccccccc]}
(0,0)& \tilde{\mu} & \tilde{\lambda} &  &  &  &  &  &  &  & \\ 
(1,0)&  &  & \tilde{\lambda} & \tilde{\mu} &  &  &  &  &  & \\ 
(2,0)&  &  & \tilde{\lambda} &  & \tilde{\mu} &  &  &  &  & \\ 
(0,1)&  &  &  & \tilde{\mu} & \tilde{\lambda} &  &  &  &  & \\ 
(1,1)&  &  &  &  &  & \tilde{\lambda} & \tilde{\mu} &  &  & \\ 
(2,1)&  &  &  &  &  & \tilde{\lambda} &  & \tilde{\mu} &  & \\ 
(0,2)&  &  &  &  &  &  & \tilde{\mu} & \tilde{\lambda} &  & \\ 
(1,2)&  &  &  &  &  &  &  &  & \tilde{\lambda} & \tilde{\mu} \\ 
(2,2)&  &  &  &  &  &  &  &  & \tilde{\lambda} & \tilde{\mu} \\
(0,3)&  &  &  &  &  &  &  &  &  & \tilde{\lambda}+\tilde{\mu} \\
\end{block}
\end{blockarray} 
\end{equation*}

Assuming we start in state $(0,0)$, the starting vector is $s=[1,0,0,0,0,0,0,0,0,0]$.
By conditioning on \numberevents$=2$, we obtain the distribution over the state space given by:
\begin{align*}
    s\bar{P}^2&=[\mu^2,\lambda\mu,\lambda^2,\lambda\mu,0,0,0,0,0,0]
\end{align*}
Since the first three states ($\{(0,0),(1,0),(2,0)\}$) correspond to $\textrm{\production{t+L}} - \textrm{\production{t}}=0$, the following three states ($\{(0,1),(1,1),(2,1)\}$) to $\textrm{\production{t+L}} - \textrm{\production{t}}=1$, etc. we obtain the following conditional probabilities:
\begin{align*}
    P(\textrm{\production{t+L}} - \textrm{\production{t}}=0|I_1^P(t)=0,\textrm{\numberevents}=2)&=\mu^2+\lambda\mu+\lambda^2\\
    P(\textrm{\production{t+L}} - \textrm{\production{t}}=1|I_1^P(t)=0,\textrm{\numberevents}=2)&=\lambda\mu\\
    P(\textrm{\production{t+L}} - \textrm{\production{t}}=2|I_1^P(t)=0,\textrm{\numberevents}=2)&=0
\end{align*}
Similarly, we can obtain the conditional probability distribution over the production of \compone{} for different values of \numberevents. Using Equation \eqref{eq: unconditional_distribution} we can then determine the distribution of the production given that there are $I_1^P(t)=0$ customer orders waiting in the queue.

When starting with $i=1$ and starting with $i=2$ customer orders in the queue, we can repeat the calculations by using a different vector $s$ (i.e., $s=[0,1,0,0,0,0,0,0,0,0]$ and $s=[0,0,1,0,0,0,0,0,0,0]$, respectively).

\subsection{Costs Calculations}
Now that we have formulated the inventory policy with target inventory position $\tilde{IP}_2(i)$ when $I_1^P(t)=i$ and modeled the production of module \compone, we want to compare the costs of this policy to the costs of a policy with a fixed target inventory position. First, we formulate the expected costs of the proposed policy in Equation \eqref{eq: expected costs state-dependent policy}. By conditioning on the value of $I_1^P(t)=i$ we determine the expected costs for every $i$. Since the stationary probability of having $i$ customers waiting in an M/M/1 queue is equal to $P(I_1^P(t)=i)=(1-\rho)\rho^i$, we can then calculate the overall expected costs.
\begin{align}\label{eq: expected costs state-dependent policy}
    E[\tilde{C}(t+L)|\tilde{IP}_2(t)]&= \sum_{i=0}^{\infty} P(I_1^P(t)=i) E[\tilde{C}(t+L)|\tilde{IP}_2(t), I_1^P(t)=i]\\
    &= \sum_{i=0}^{\infty} (1-\rho)\rho^i \sum_{j=0}^{X^U} P(\textrm{\production{t+L}} - \textrm{\production{t}}=j|I_1^P(t)=i) \left( (b+h_1)(j - \tilde{IP}_2(i) )^+ + h_2 (\tilde{IP}_2(i) - j)^+ \right)
\end{align}
For an inventory policy with fixed target inventory position $IP_2$, the expected costs are given in Equation \eqref{eq: expected costs fixed policy}. Without taking into account the current state of the queue, the output process of the M/M/1 queue during $[t,t+L]$ is a Poisson process with rate $\lambda L$. Therefore,  $P(\textrm{\production{t+L}} - \textrm{\production{t}}=j)=\frac{(\lambda L)^j}{j!}e^{-\lambda L}$.
\begin{align}\label{eq: expected costs fixed policy}
    E[\tilde{C}(t+L)|IP_2(t)]= \sum_{j=0}^{X^U} P(\textrm{\production{t+L}} - \textrm{\production{t}}=j) \left( (b+h_1)(j - IP_2 )^+ + h_2 (IP_2 - j)^+ \right)
\end{align}

\subsection{Computational Analysis}\label{sec: computational analysis}
In this section, we evaluate the effectiveness of the inventory policy that takes into account the customers currently waiting at the in-house production facility of module \compone{} by means of a computational experiment. The results of this analysis are provided in \S\ref{sec: policy evaluation}. First, in \S\ref{sec: example policy} we illustrate how the procedure described in \S \ref{sec: model} results in a table of target IPs for different numbers of waiting customer orders.

\subsubsection{Illustrative Example of Inventory Policy}\label{sec: example policy}
Consider the following numerical example: $\lambda=0.8$, $\mu=1$, $L=4$, $h_1=4$, $h_2=1$, $b=5$, giving critical fractile $CF=0.9$.
In Table \ref{tab: target_IP} we tabulate the target IP for module \comptwo{} for different values of $I_1^P(t)$.
\begin{table}[ht]
    \centering
    \caption{Target inventory position for different values of $I_1^P(t)$}
    \begin{tabular}{|r|r|}
        \hline
        $I_1^P(t)$ & Target IP\\
        \hline
        0	&	4	\\
        1	&	4	\\
        2   &   5   \\
        3   &   6   \\
        4   &   6   \\
        5   &   6   \\
        
        $\geq 6$	&	7	\\
        \hline
    \end{tabular}
    \label{tab: target_IP}
\end{table}

Without taking into account the current state of the queue, the output process of the M/M/1 queue during $[t,t+L]$ is a Poisson process with rate $\lambda L$. Therefore,  $P(\textrm{\production{t+L}} - \textrm{\production{t}}=j)=\frac{(\lambda L)^j}{j!}e^{-\lambda L}$. Using this information, we can determine the target inventory position for which the cumulative probability reaches the critical fractile. In the given example, this leads to a target inventory position of 6.


\subsubsection{Policy Evaluation}\label{sec: policy evaluation}


In order to assess the value of taking into account information on outstanding orders in the order policy and its sensitivity to various model parameters, we perform a full factorial experiment. In our experiment, we vary the rate of production relative to the rate of incoming orders, the lead time of module \comptwo, the holding costs of module \compone{} relative to those of module \comptwo{} and the waiting time costs for customers. The setup of the experiment is given in Table \ref{tab: experiment setup}. We set $\mu=1$ and $h_2=1$. In total we have $2\cdot3^3=54$ instances. 

\begin{table}[ht]
    \centering
    \caption{Parameter settings for experiments}
    \begin{tabular}{ll}
        \hline
        Parameter & Values  \\
        \hline
        $\lambda$ & 0.8, 0.9\\
        $L$ & 2, 3, 4\\
        $h_1$ & 1, 2, 4\\
        $b$ & 1, 5, 10\\
        \hline
    \end{tabular}
    \label{tab: experiment setup}
\end{table}

For each instance we calculate the expected costs at time $t+L$ under the reorder policy with fixed target inventory position and the policy with state-dependent target inventory position and the cost savings that can be achieved by taking into account information on customer orders. The parameters given in Table \ref{tab: experiment setup} are fixed one at a time and the remaining parameters are varied. The summary statistics are provided in Table \ref{tab: experiment results}. 
The last row provides the summary over all 54 instances.

\begin{table}[ht]
    \centering
    \caption{Results of full factorial experiment continuous time model}
    \begin{tabular}{lrrrrrr}
    \hline
    & & \multicolumn{1}{p{3cm}}{\centering Average cost\\ fixed IP}   & \multicolumn{1}{p{3cm}}{\centering Average cost\\ state-dep. IP} & Avg saving (\%) & Max saving (\%)  & Min saving (\%) \\
    \hline
    $\lambda$	&	0.8	&	2.71	&	2.30	&	14.65	&	18.27	&	8.72	\\
	&	0.9	&	2.87	&	2.31	&	19.34	&	22.07	&	15.02	\\ \hline
$L$	&	2	&	2.34	&	1.95	&	16.05	&	20.90	&	8.72	\\
	&	3	&	2.82	&	2.33	&	16.95	&	20.77	&	10.68	\\
	&	4	&	3.22	&	2.63	&	17.98	&	22.07	&	13.33	\\ \hline
$h_1$	&	1	&	2.57	&	2.14	&	16.08	&	21.23	&	8.72	\\
	&	2	&	2.77	&	2.28	&	17.45	&	22.07	&	11.62	\\
	&	4	&	3.03	&	2.49	&	17.46	&	21.85	&	12.57	\\ \hline
$b$	&	1	&	2.14	&	1.79	&	15.68	&	20.85	&	8.72	\\
	&	5	&	2.88	&	2.39	&	16.95	&	20.90	&	10.51	\\
	&	10	&	3.35	&	2.74	&	18.36	&	22.07	&	15.80	\\ \hline
All	&		&	2.79	&	2.31	&	17.00	&	22.07	&	8.72	\\ \hline
    \end{tabular}
    \label{tab: experiment results}
\end{table}

According to the results provided in Table \ref{tab: experiment results}, on average 17.00\% cost reduction can be achieved by letting the target IP depend on the number of waiting customers. The minimum and maximum cost reduction are 8.72\% and 22.07\%, respectively. 
The benefit of using the available information on the number of waiting customer orders increases considerably in the arrival rate. When the arrival rate is low compared to the service rate, the number of waiting customers is likely to be small and the Poisson process with rate $\lambda L$ will be a good approximation of the production of module \compone, hence there is less value in using this information. However, when $\lambda$ is close to $\mu$, the queue of waiting customers may be more substantial and the information on the number of outstanding customer orders becomes more useful.
Furthermore, the value of this information is higher when the holding cost of module \compone{} and/or the customer waiting costs are high relative to the holding costs of module \comptwo{} or when the lead-time of module \comptwo{} is high. Overall, Table \ref{tab: experiment results} shows that including the information on waiting customer orders in the order policy leads to substantial savings.

\section{Discrete Time Model}\label{sec: discrete model}
In the previous sections we have shown that an inventory policy with a state-dependent target inventory position is optimal and can generate considerable savings compared to a policy with a fixed target inventory position in a continuous time model. Such a model may not be suitable for all assembly systems, for example when using periodic review. Therefore, we  will now consider a model in discrete time.
We will again consider the assembly system shown in Figure \ref{fig: sketch assembly system}, consisting of two modules that need to be merged in a single end-product. Module \comptwo{} is sourced from a supplier with a lead-time of $L$ periods. For the other module there is in every period available production capacity that may either be fixed or random. 
We denote the (random) number of units that can be produced per period by $C$ and demand per period by $D$. 

The production of module \compone{} can still be modeled as a Markov process with states $(i,j)$, where $i$ denotes the number of outstanding customer orders and $j$ the number of units produced. Every period production is equal to the minimum of available capacity ($C$) and available customer orders  consisting of both new demand and outstanding orders ($i+D$). If the total number of orders exceeds available capacity, the remaining orders will still be outstanding orders at the beginning of the next period. This means that we have the following transition:
\begin{align*}
    (i,j) & \rightarrow \left((i+D-C)^+,j+\min\{i+D,C\}\right).
\end{align*}
The transition probabilities can be determined based on the distributions of demand and capacity. 

The cost structure remains the same as in the continuous time model. In every period holding costs $h_1$ are incurred for finished modules \compone{} that need to be stored and similarly we have holding costs $h_2$ for \comptwo. Additionally, there is a per-period back-order cost $b$. Costs per period are as given in Equation \eqref{eq: costs}. Using similar reasoning as for the continuous time model, we can write the myopic inventory policy as:
\begin{align}\label{eq: minimization periodic}
    \min_{IP_2(t)} E[(b+h_1)(\textrm{\production{t+L}} - \textrm{\production{t}} - IP_2(t) )^+ + h_2 (IP_2(t) - (\textrm{\production{t+L}} - \textrm{\production{t}}))^+|I_1^P(t)].
\end{align}

We can again show that the myopic inventory policy, in which the target inventory position is determined for every period separately, is optimal. For this we prove Theorems \ref{thm: periodic monotonicity} and \ref{thm: periodic max_increase}, which are similar to Theorems \ref{thm: monotonicity} and \ref{thm: max_increase} in the continuous model. These theorems show that the target inventory position of module \comptwo{} is non-decreasing in the number of outstanding customer orders and increases by at most 1 as the number of outstanding orders increases by 1.

\begin{theorem}\label{thm: periodic monotonicity}
When both demand $D$ and capacity $C$ per period are uncertain, $\tilde{IP}_2(I_1^P(t))$ is monotonically non-decreasing in $I_1^P(t)$.
\end{theorem}

\begin{theorem}\label{thm: periodic max_increase}
When both demand $D$ and capacity $C$ per period are uncertain, an additional customer in the system increases the target inventory position of module \comptwo{} by at most 1.
\end{theorem}

Similar to the continuous time model, we can use Theorems \ref{thm: periodic monotonicity} and \ref{thm: periodic max_increase} to conclude that the proposed policy is well-defined also for the discrete time model. Furthermore, we can use the same reasoning as in \S\ref{sec: inventory policy c2} to conclude that the myopic policy for module \comptwo{} is optimal in discrete time in Theorem \ref{thm: periodic optimality}.

\begin{theorem}\label{thm: periodic optimality}
When both demand $D$ and capacity $C$ per period are uncertain, the myopic inventory policy for module \comptwo{} is optimal.
\end{theorem}

In summary, the proposed inventory policy is not only optimal in the continuous time $M/M/1$ model, but all analytical results continue to hold when we consider a setting with discrete time periods and a capacity per period that is either fixed or random.

\subsection{Computational Analysis}
Now that we have established that the proposed myopic policy is optimal also in the discrete time case, we will evaluate the effectiveness of this policy compared to a base-stock policy with a fixed base-stock level. In our experiment, we vary the lead time of module \comptwo, the holding costs of module \compone{} relative to those of module \comptwo{} and the waiting time costs for customers. We set $h_2=1$ and use the same parameter values for $L$, $h_1$ and $b$ as given in Table \ref{tab: experiment setup}. Furthermore, we consider different cases with respect to the distributions for demand and capacity:

\begin{table}[ht]
    \centering
    \begin{tabular}{ccc}
        Case 1: &  $P(D=d)=\left\{\begin{matrix}
0.4 & \text{if }d=1\\ 
0.4 & \text{if }d=2\\ 
0.2 & \text{if }d=3
\end{matrix}\right.$, &  $P(C=c)=\left\{\begin{matrix}
0.7 & \text{if }c=2\\ 
0.3 & \text{if }c=3
\end{matrix}\right.$\\
& & \\
        Case 2: &  $P(D=d)=\left\{\begin{matrix}
0.1 & \text{if }d=0\\
0.2 & \text{if }d=1\\ 
0.2 & \text{if }d=2\\ 
0.25 & \text{if }d=3\\ 
0.15 & \text{if }d=4\\ 
0.1 & \text{if }d=5
\end{matrix}\right.$, &  $P(C=c)=\left\{\begin{matrix}
0.2 & \text{if }c=2\\ 
0.45 & \text{if }c=3\\
0.35 & \text{if }c=4
\end{matrix}\right.$
    \end{tabular}
\end{table}

In total we thus have again $2\cdot 3^3=54$ instances. For each instance we calculate the expected costs at time $t+L$ under both the reorder policy with fixed target inventory position and the policy with state-dependent target inventory position. Subsequently, we determine the cost savings achieved by using the state-dependent policy. We again perform a full factorial experiment, where we vary the parameter values or demand and capacity distribution cases one by one. 
We determine the distribution of the number of outstanding orders for module \compone{} by starting with an empty system and taking the distribution over the number of outstanding orders after $N$ periods, where $N$ is selected such that the probabilities have stabilized.
The results of the experiment are given in Table \ref{tab: discrete experiment results}. 

Similar to the continuous time model, we observe an increase in the average savings when considering the information on outstanding orders of module \compone{} in your inventory policy for module \comptwo{} as the holding costs $h_1$ increase. There are also some differences compared to the results for the continuous time model. In the continuous time case, the average savings percentage increased as the lead-time increased, whereas in the discrete time case we observe the opposite effect. When $L=2$ the average savings from incorporating the information on outstanding orders for module \compone{} are 11.70\% and for $L=4$ this has decreased to 9.05\%. This is due to the fact that we now have discrete distributions for demand and capacity with a finite number of outcomes. Therefore, there is no possibility of high peaks in the number of demands occurring during a certain time period and, hence, in the discrete case, the value of synchronization reduces as the lead-time grows large.
Furthermore, we no longer observe an increasing trend in the average savings percentage as customer waiting costs $b$ increase. When we consider the two cases with respect to the distributions of demand and capacity, we observe that the expected demand relative to the expected capacity is comparable in both cases. In the second case the variation in both capacity and demand per period is larger. This has a large effect on the average savings, as in case 1 the average savings are equal to 6.84\% and in case 2 the average savings are 13.47\%. 

\begin{table}[ht]
    \centering
    \caption{Results of full factorial experiment discrete time model}
    \begin{tabular}{lrrrrrr}
    \hline
    & & \multicolumn{1}{p{3cm}}{\centering Average cost\\ fixed IP}   & \multicolumn{1}{p{3cm}}{\centering Average cost\\ state-dep. IP} & Avg saving (\%) & Max saving (\%)  & Min saving (\%) \\
    \hline
Case	&	1	&	1.55	&	1.45	&	6.84	&	25.02	&	0.02	\\
	&	2	&	2.85	&	2.45	&	13.47	&	26.28	&	4.24	\\ \hline
$L$	&	2	&	1.78	&	1.53	&	11.70	&	26.28	&	0.02	\\
	&	3	&	2.21	&	1.96	&	9.71	&	21.28	&	1.21	\\
	&	4	&	2.62	&	2.35	&	9.05	&	18.79	&	1.77	\\ \hline
$h_1$	&	1	&	2.06	&	1.83	&	9.50	&	21.94	&	0.02	\\
	&	2	&	2.14	&	1.93	&	8.53	&	20.78	&	0.04	\\
	&	4	&	2.40	&	2.08	&	12.42	&	26.28	&	0.19	\\ \hline
$b$	&	1	&	1.76	&	1.59	&	9.44	&	25.02	&	2.60	\\
	&	5	&	2.28	&	2.00	&	11.71	&	26.28	&	1.36	\\
	&	10	&	2.56	&	2.26	&	9.31	&	21.94	&	0.02	\\ \hline
All	&		&	2.20	&	1.95	&	10.15	&	26.28	&	0.02	\\ \hline
    \end{tabular}
    \label{tab: discrete experiment results}
\end{table}

\section{Discussion and Conclusion}\label{sec: conclusions}
We have examined an assembly problem where an OEM builds a final product consisting of two modules, one of which is made-to-order by the OEM itself. Since production of the make-to-order module commences as soon as a customer order arrives and the final product can be assembled as soon as both modules are available, intermediate stocks and thus costs are controlled by the order policy of the module sourced from the supplier. 
Since the number of orders waiting for production of the in-house produced module gives an indication of the demand for the other module one lead-time from now, this information can be useful in determining the inventory policy for this module. Therefore, we consider an inventory policy where the target inventory position depends on the state of the production system for the in-house produced module. When there is a large number of orders waiting for production of the module, production during the lead-time of the other module is likely to be higher than in case none or very few orders are waiting. Consequently, the target inventory position of the supplier-sourced module will be higher when there are many customer orders waiting in the in-house production system.

We show that under this policy the target inventory position is monotonically increasing in the number of customer orders in the queue. 
Additionally, we show optimality of this policy both in continuous and discrete time.
In support of these analytical results and to illustrate the proposed policy, we conducted a computational analysis. In this analysis we performed a full factorial experiment to evaluate the benefit of taking information on outstanding customer orders into account. We show that using this information can lead to considerable savings. Furthermore, we assess the sensitivity to various model parameters and show that especially the arrival rate of customer orders relative to the production rate has a large influence on the savings. 


\section*{Acknowledgements}
This work is part of the research programme Complexity in high-tech manufacturing
with project number 439.16.121, which is (partly) financed by the Dutch Research Council (NWO).

\printbibliography

\begin{appendices}\label{sec: appendix}
\section{Proofs}

\paragraph{Proof of Theorem \ref{thm: monotonicity}}
\begin{proof}
By conditioning on a specific sequence of events, i.e. customer arrivals and finished production of a module, we show that it holds for any possible event sequence and thus also for the expectation over all possible event sequences.

Denote the following:
\begin{align*}
    J(i) &= \textrm{\# finished modules \compone{} when starting at } (i,0)\\
    J(i+1) &= \textrm{\# finished modules \compone{} when starting at } (i+1,0)
\end{align*}

We compare $J(i)$ and $J(i+1)$, so we consider starting at $(i,0)$ and $(i+1,0)$. There are two transitions, namely arrival of a new customer order to the queue and finished production of a module. This gives the following cases to consider.
\begin{enumerate}
    \item \textbf{Arrival of customer order} gives transitions $(i,0)\rightarrow(i+1,0)$ and $(i+1,0)\rightarrow(i+2,0)$. 
    \item \textbf{Finished production \& $i>0$} gives transitions $(i,0)\rightarrow(i-1,1)$ and $(i+1,0)\rightarrow(i,1)$.
    \item \textbf{Finished production \& $i=0$} gives transitions $(i,0)=(0,0)\rightarrow(0,0)$ and $(i+1,0)=(1,0)\rightarrow(0,1)$.
\end{enumerate}
For cases 1 and 2, the difference in number of customers in the queue remains 1 and the difference between the number of finished modules remains 0. Thus, $J(i)$ and $J(i+1)$ stay equal in these two cases. For case 3, the number of customers in the queue becomes equal, but the difference in the number of finished modules increases by one and will thus always stay one ahead from now on.

After these transitions, the next transition also falls within one of the three cases discussed above. Since $J(i+1)$ either remains equal to $J(i)$ or stays one finished module ahead of $J(i)$, it follows that $J(i)\leq J(i+1)$. This means that for every number of finished modules $j$, $P(J(i+1)>j)\geq P(J(i)>j)$. From this we can conclude that $J(i+1)$ is statistically larger than $J(i)$ and thus that the target inventory position of module \comptwo{}  when $i+1$ customer orders are waiting in the queue is equal to or larger than in case $i$ customer orders are waiting in the queue.
\end{proof}

\paragraph{Proof of Theorem \ref{thm: max_increase}}
\begin{proof}
This can be proven along the same lines as Theorem \ref{thm: monotonicity}. Using the definitions of $J(i)$ and $J(i+1)$, it follows that $J(i)+1$ equals the \# finished modules \compone{} when starting at $(i,1)$.
We will consider the case of $i$ orders in the queue and already one finished module and the case with $i+1$ customers in the queue and no finished modules. This means that we compare $J(i)+1$ and $J(i+1)$.

Again, there are two transitions, namely arrival of a new customer to the queue and departure of a customer from the system, corresponding to finished production of a module. This gives the following cases to consider.
\begin{enumerate}
    \item \textbf{Arrival} gives transitions $(i,1)\rightarrow(i+1,1)$ and $(i+1,0)\rightarrow(i+2,0)$. 
    \item \textbf{Departure \& $i>0$} gives transitions $(i,1)\rightarrow(i-1,2)$ and $(i+1,0)\rightarrow(i,1)$.
    \item \textbf{Departure \& $i=0$} gives transitions $(i,1)=(0,1)\rightarrow(0,1)$ and $(i+1,0)=(1,0)\rightarrow(0,1)$.
\end{enumerate}
For cases 1 and 2, the difference in number of customers in the queue remains 1 and the difference between the number of exits remains 1. Thus, $J(i)+1$ stays one finished module ahead of $J(i+1)$ in these two cases. In case 3, both the number of customers in the queue and the number of finished modules after the transition are equal, meaning that $J(i)+1$ and $J(i+1)$ are equal and will follow the same trajectory when facing new arrivals or departures. 

After these transitions, the next transition also falls within one of the three cases discussed above. Since $J(i)+1$ either stays one finished module ahead of $J(i+1)$ or they become equal and follow the same trajectory, it follows that $J(i+1)\leq J(i)+1$. This means that for every number of exits $j$, $P(J(i+1)>j)\leq P(J(i)+1>j)$ and thus that $J(i+1)$ is statistically smaller than $J(i)+1$. 

From this we can conclude that the number of finished products during $[t,t+L]$ when starting with $i+1$ orders in the queue is at most one larger than in case there are $i$ orders in the queue. Therefore, the target inventory position of module \comptwo{} is also at most one higher.
\end{proof}

\paragraph{Proof of Theorem \ref{thm: optimality}}
\begin{proof}
The optimal inventory policy, denoted by $\pi$, for module \comptwo{} is the one that minimizes total costs:
\begin{equation*}
    \min_{\pi}\int_0^{T} \mathbb{E}[\tilde{C}(t)|\pi]dt
\end{equation*}
where $T$ denotes the planning horizon.

Given the lead-time $L$, the inventory levels and thus costs at time $t$ are affected by the decisions made at time $t-L$ and corresponding inventory position. Since every inventory policy for module \comptwo{} has a corresponding inventory position, it holds that
\begin{equation*}
    \min_{\pi}\int_0^{T} \mathbb{E}[\tilde{C}(t)|\pi]dt \geq \int_0^{T} \min_{IP_{t-L}} \mathbb{E}[\tilde{C}(t)|IP_{t-L}]dt
\end{equation*}
Hence, a lower bound on the costs is obtained for the policy in which the  optimal inventory position is selected at every time $t$. This lower bound is attained by the proposed myopic inventory policy for module \comptwo. Therefore, it can be concluded that the proposed policy is optimal.
\end{proof}

\paragraph{Proof of Lemma \ref{lem: upper bound}}
\begin{proof}
Cantelli's inequality states that for $X$ with mean $\hat{\mu}$ and variance $\hat{\sigma}^2$, $P(X\geq r)\leq \frac{\hat{\sigma}^2}{\hat{\sigma}^2+(r-\hat{\mu})^2}$ for $r>\hat{\mu}$. Since we have $X\sim Poisson((\lambda+\mu)L)$, $\hat{\mu}=(\lambda+\mu)L$ and $\hat{\sigma}^2=(\lambda+\mu)L$, this gives $P(X\geq r)\leq \frac{(\lambda+\mu)L}{(\lambda+\mu)L+(r-(\lambda+\mu)L)^2}$. To find $r$ such that $P(X\geq r)$, we need to solve $\frac{(\lambda+\mu)L}{(\lambda+\mu)L+(r-(\lambda+\mu)L)^2}=\epsilon$. This yields $r=(\lambda+\mu)L+\sqrt{\left(1-\frac{1}{\epsilon}\right)(\lambda+\mu)L}$. Therefore, $P(X\ge X^U)\leq \epsilon$ holds if $X^U=(\lambda+\mu)L+\sqrt{\left(1-\frac{1}{\epsilon}\right)(\lambda+\mu)L}$.
\end{proof}

\paragraph{Proof of Theorem \ref{thm: periodic monotonicity}}
\begin{proof}
The proof is along the same lines as for Theorem \ref{thm: monotonicity}. We again compare $J(i)$ and $J(i+1)$. When we are in state $(i,j)$, meaning that there are $i$ product orders in the system and $j$ units produced so far, and demand $D$ occurs, there are $i+D$ units to be produced. Since the available capacity is $C$ units, production is equal to $\min\{i+D,C\}$ and the remaining number of back-orders equals $(i+D-C)^+$. 

By conditioning on demand $D=d$ and capacity $C=c$, the following transition occurs:
\begin{align*}
    (i,j) & \rightarrow \left((i+d-c)^+,j+\min\{i+d,c\}\right).
\end{align*}
We consider the following cases:
\begin{enumerate}
\item $i+d<c$ (so $i+1+d\leq c$): in this case we have transitions
\begin{align*}
    (i,j) & \rightarrow (0,j+i+d) \textrm{ and }\\
    (i+1,j) & \rightarrow (0,j+i+d+1)
\end{align*}
Since all product orders can be satisfied, the number of outstanding product orders reduces to zero for both $(i,j)$ and $(i+1,j)$. Therefore, from now on production in any period will be the same and $J(i+1)$ will stay one unit ahead of $J(i)$.
\item $i+d\geq c$: gives transitions
\begin{align*}
    (i,j) & \rightarrow (i+d-c,j+c) \textrm{ and }\\
    (i+1,j) & \rightarrow (i+1+d-c,j+c)
\end{align*}
In this case the difference in the number of outstanding order stays 1 and the number of produced items remains equal. Hence, the following period starts of with the same situation as the current period.
\end{enumerate}
Since $J(i+1)$ either stays equal to $J(i)$ or stays one finished module ahead of $J(i)$, it follows that $J(i)\leq J(i+1)$. This means that for every number of finished modules $j$, $P(J(i+1)>j)\geq P(J(i)>j)$. From this we can conclude that $J(i+1)$ is statistically larger than $J(i)$ and thus that the target inventory position of module \comptwo{}  when $i+1$ customer orders are in the system is equal to or larger than in case $i$ customer orders are in the system.
\end{proof}

\paragraph{Proof of Theorem \ref{thm: periodic max_increase}}
\begin{proof}
The proof is along the same lines as for Theorem \ref{thm: monotonicity} and \ref{thm: periodic monotonicity}. We again compare $J(i+1)$ and $J(i)+1$. 

We again condition on demand $D=d$ and capacity $C=c$ and consider the following cases:
\begin{enumerate}
\item $i+d<c$ (so $i+1+d\leq c$): in this case we have transitions
\begin{align*}
    (i,j+1) & \rightarrow (0,j+1+i+d) \textrm{ and }\\
    (i+1,j) & \rightarrow (0,j+i+d+1)
\end{align*}
The number of outstanding product orders again reduces to zero for both $(i,j+1)$ and $(i+1,j)$. Also, the number of produced items becomes equal. Therefore, from now on production in any period will be the same and $J(i+1)$ will remain the same $J(i)$.
\item $i+d\geq c$: gives transitions
\begin{align*}
    (i,j+1) & \rightarrow (i+d-c,j+1+c) \textrm{ and }\\
    (i+1,j) & \rightarrow (i+1+d-c,j+c)
\end{align*}
In this case the difference in the number of outstanding order stays and the difference in the number of produced items both remain equal. Hence, the following period starts of with the same situation as the current period.
\end{enumerate}
After these transitions, the next transition also falls within one of the two cases discussed above. Since $J(i)+1$ either stays one finished module ahead of $J(i+1)$ or they become equal and follow the same trajectory, it follows that $J(i+1)\leq J(i)+1$. This means that for every number of exits $j$, $P(J(i+1)>j)\leq P(J(i)+1>j)$ and thus that $J(i+1)$ is statistically smaller than $J(i)+1$. 

From this we can conclude that the number of finished products when starting with $i+1$ orders in the queue is at most one larger than in case there are $i$ orders in the queue. Therefore, the target inventory position of module \comptwo{} is also at most one higher.
\end{proof}

\paragraph{Proof of Theorem \ref{thm: periodic optimality}}
\begin{proof}
The optimal inventory policy, denoted by $\pi$, for module \comptwo{} is the one that minimizes total costs:
\begin{equation*}
    \min_{\pi}\sum_0^{T} \mathbb{E}[\tilde{C}(t)|\pi]dt
\end{equation*}
where $T$ denotes the planning horizon.

Given the lead-time of $L$ periods, the inventory levels and thus costs in period $t$ are affected by the decisions made in period $t-L$ and corresponding inventory position. Since every inventory policy for module \comptwo{} has a corresponding inventory position, it holds that
\begin{equation*}
    \min_{\pi}\sum_0^{T} \mathbb{E}[\tilde{C}(t)|\pi]dt \geq \sum_0^{T} \min_{IP_{t-L}} \mathbb{E}[\tilde{C}(t)|IP_{t-L}]dt
\end{equation*}
Hence, a lower bound on the costs is obtained for the policy in which the  optimal inventory position is selected in every period $t$. This lower bound is attained by the proposed myopic inventory policy for module \comptwo. Therefore, it can be concluded that the proposed policy is optimal.
\end{proof}

\end{appendices}
\end{document}